%% file: main.tex
\documentclass[letterpaper,10 pt,twoside]{ieeeconf}

\usepackage{amssymb}
\usepackage{cite}
\usepackage{layouts}
\usepackage{soul}

\IEEEoverridecommandlockouts

\def\BibTeX{{\rm B\kern-.05em{\sc i\kern-.025em b}\kern-.08em
    T\kern-.1667em\lower.7ex\hbox{E}\kern-.125emX}
}%

\input{ls/generic_article.tex}
\input{ls/pm_specific_setup.tex}

\title{Persistent Monitoring Trajectory Optimization in Partitioned Environments}

\author{\IEEEauthorblockN{}
\IEEEauthorblockA{\textit{Systems Engineering} \\
\textit{Boston University}\\
Massachusetts, USA
}
}

\author{Jonas Hall$^{1}$, Christos G. Cassandras$^{1,2}$, Sean B. Andersson$^{1,3}$
\thanks{This work was supported in part by the NSF through ECCS-1931600}%
\thanks{$^{1}$Division of Systems Engineering, Boston University, USA}%
\thanks{$^{2}$Department of Electrical and Computer Engineering}%
\thanks{$^{3}$Division of Mechanical Engineering, Boston University, USA}%
\thanks{{\tt\small \{hallj, cgc, sanderss\}@bu.edu}}
}%

\newcommand{\scalefactor}{0.810}


\begin{document}

\clearpage\maketitle

\begin{abstract}
    We consider the problem of using an autonomous agent to persistently monitor a collection of targets distributed in a given environment. We generalize existing work by allowing the agent's dynamics to vary throughout the environment, leading to a hybrid dynamical system. This introduces an additional layer of complexity towards the planning portion of the problem: we must not only identify in which order to visit the points of interest, but also in which order to traverse the regions. We propose a tailored global path planner and prove that it is not only probabilistically complete, but converges in probability to a time-optimal solution. We then design an offline sequence planner together with an online trajectory optimizer. Simulations validate the results.
\end{abstract}

\section{Introduction}

\ac{pm} describes a broad class of problems in which an agent moves through an environment to collect information about specific targets over time. It is applicable across a wide range of applications, such as ocean monitoring~\cite{smith2011persistentOcean}, forest fire surveillance~\cite{casbeer2006cooperative}, tracking of individual biological macromolecules~\cite{pinto2021tracking}, and data harvesting~\cite{zhu2022control}. As a specific motivating application, consider a disaster scenario. An efficient response to a catastrophic event such as an earthquake, hurricane, or tsunami requires the persistent and simultaneous state estimation of many locations in order to make time-sensitive resource distribution decisions. The affected area may consist of various types of terrain, each having its own characteristics. For instance, some parts may be urban and thus contain obstacles, others may be coastal regions with the presence of strong winds. In this paper, we propose an online trajectory optimization scheme to minimize the average estimation error of the target states for piecewise continuous agent dynamics. These hybrid dynamics have the capability of capturing the various terrains the agent must move through; in turn, this requires the extension of standard trajectory optimizers for \ac{pm} in order to exploit the local structures.

In this paper, we assume that there are multiple targets to be monitored, each of which is located in a different region. We model the \textit{state} of each target with a stochastic \ac{lti} system, and assume the agent has the ability to interact with each target in its vicinity by taking noisy measurements of this state; this typical setting can also be found in \ac{pm}~\cite{lan2013planning,chen2019deep,ha2019periodic}. To obtain good estimates of the state of a target, the agent must remain nearby to collect measurements over time. However, since the agent is tasked with monitoring all states, it cannot spend too much time at each individual location since the accuracy of the estimates of the states that are not observed decays over time. Thus, we want to optimize both the time spent at the targets and the order in which to visit them.

The optimization of the visiting sequence is often approached by abstracting the problem and casting it into a \ac{tsp} or \ac{vrp} in an offline preparation phase~\cite{hari2020optimal,stump2011multi,yu2018optimal}. This requires the estimation of travel times between targets, which is itself a challenging problem in our setting of hybrid agent dynamics. It is conceivable to tackle this problem with modern optimization-based solvers~\cite{nurkanovic2022nosnoc}. However, they are prone to getting stuck at local minima and require careful initialization. In contrast, sampling-based methods often come with theoretical guarantees, such as probabilistic convergence to optimal solutions~\cite{karaman2010incremental}. Their main drawback is the computational cost. However, since the computation of the visiting sequence is already performed off-line, utilizing a sampling-based solver does not present a bottleneck at this stage, allowing us to exploit its preferred robustness properties.

Building upon our prior work~\cite{hall2023bilevel}, we optimize the agent trajectory on-line, by optimizing parameters of a decomposition describing a periodic trajectory that realizes the computed visiting sequence. The present paper differs from~\cite{hall2023bilevel} in many key aspects, including a more complex model of target uncertainty (described now in terms of variance of the state estimates), a more general environment (via the piecewise smooth dynamics), and a direct focus on minimizing the average uncertainty (rather than the period of the trajectory). The contributions of this paper are as follows:
\begin{itemize}
    \item a \ac{pm} formulation that allows for agent dynamics that vary across the environment.
    \item the \ac{rrbt} algorithm for global path planning in the introduced context, and a proof that it converges in probability to an optimal solution.
    \item an online optimization scheme to minimize the steady-state average estimation error.
\end{itemize}

The remainder of this paper is organized as follows. In \secref{sec:problem:formulation} we introduce the problem formulation. We characterize optimal periodic trajectories via a decomposition in \secref{sec:decomposition}.
\secref{sec:global:planning} introduces the \ac{rrbt} algorithm, and we state under which conditions the algorithm probabilistically converges to an optimal solution. \secref{sec:optimization} then provides the online optimization of the decomposed global cost function by alternating between updates of the local control segments and the global optimization parameters. The efficacy of the method is demonstrated in \secref{sec:numerical:results}, and we conclude the article in \secref{sec:future:work} while providing open questions for future work.

\section{Problem Formulation}\label{sec:problem:formulation}

\begin{figure*}
    \begin{subfigure}{0.5\linewidth}
        \centering
        \includegraphics[width=\scalefactor\linewidth]{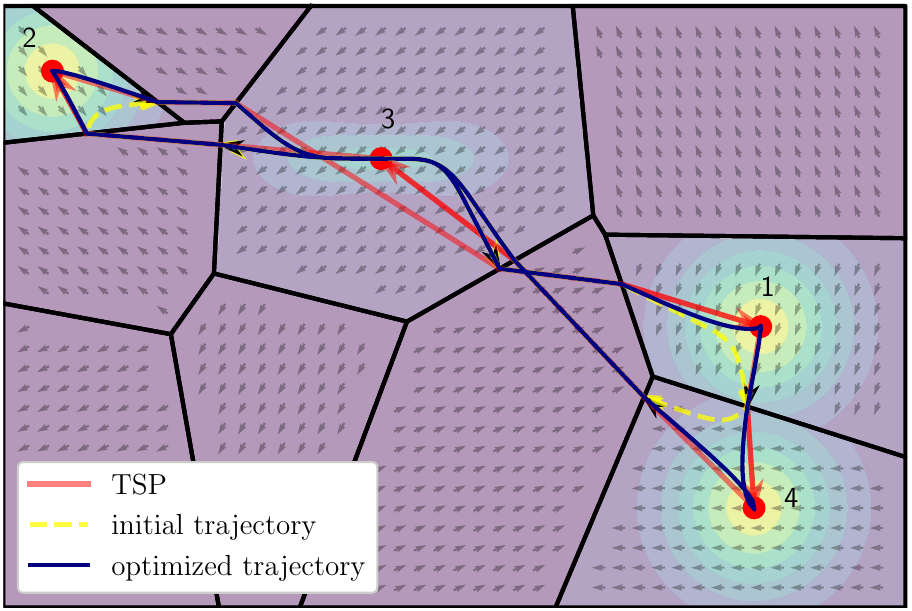}
        \caption{initial vs optimal}
        \label{fig:tsp:initial:optimal:trajectory}
    \end{subfigure}%
    \begin{subfigure}{0.5\linewidth}
        \centering
        \includegraphics[width=\scalefactor\linewidth]{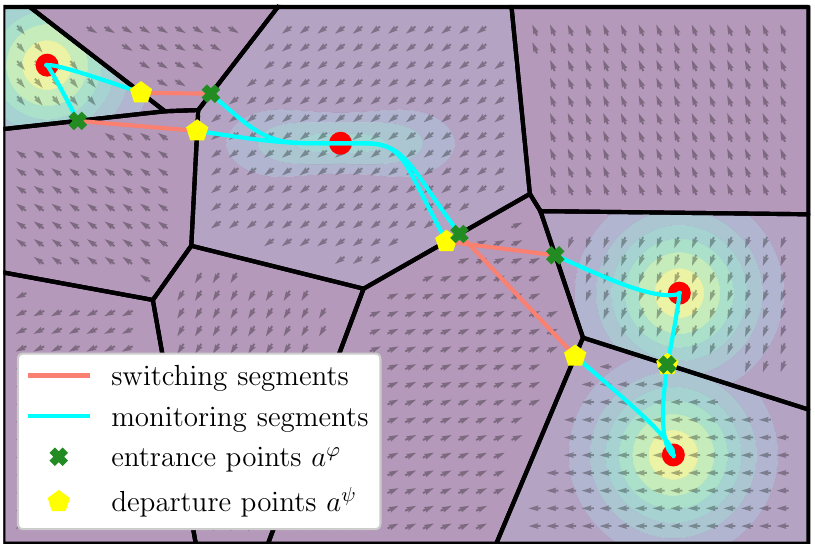}
        \caption{illustrating the decomposition}
        \label{fig:decomposition}
    \end{subfigure}
    \caption{An example of a partitioned environment with four targets and ten regions. The left plot shows the initial trajectory, the optimal trajectory, and the TSP solution. The right plot illustrates the decomposition.}
    \label{fig:example:environment}
\end{figure*}

As depicted in \figref{fig:tsp:initial:optimal:trajectory}, we consider a partitioned environment in which a single agent moves to carry out its \ac{pm} mission. In this context, the partition is driven by the agent's dynamics: these dynamics are smooth within any given region but may change abruptly when moving between regions. In the motivating disaster response example, the partition could be determined by the agent's mode of transportation, e.g., a drone that can fly over water but utilizes an auxiliary vehicle to traverse land. The partition could also be obtained as a complex model simplification into a finite collection of domains.

The agent is tasked to monitor multiple targets that are distributed across the mission space, each of which is assumed to be associated with a set of states that evolve according to stochastic \ac{lti} dynamics. The agent is equipped with a sensor that can take noisy measurements of a target's internal states. We make the simplifying assumption that each region contains at most one target, and that the agent can only sense the target when it is within the region. Let us now formalize this setup.


Denote the compact and connected mission space by $S \subset \R^2$, and assume it is partitioned into $P$ polyhedra
\begin{equation}\label{eq:partition:set}
    \partition_i = \{ x \in \R^2 \mid \partitionfunction_j\super{\top}x \leq b_j \fort j = 1, 2, \dots, P_i\},
\end{equation}
for $i = 1,2,\dots,P$, where $P_i$ denotes the number of inequalities required to describe the set $\partition_i$.

Let $\agent(t) \in \R^2$ denote the agent's position at time $t$. To capture the heterogeneous dynamics of the agent, let us define the piecewise smooth vector field
\begin{equation}\label{eq:piecewise:smooth:vf}
    f(x) = f_i(x) \ift x \in\ \interior{\partition}_i,
\end{equation}
where $\interior{\partition}_i$ denotes the interior of the region $\partition_i$. We make the assumption that the agent never moves along the boundary of a region $\partial S_i$, however, it can cross it to enter a neighboring region. This assumption is reasonable since $\cup_i \partial S_i$ has zero measure in $\S$, and any perturbation would move the agent away from the boundary. We also assume that the agent can apply a control input $u(t) \in \bar{B}_1(0)$, where $\bar{B}_1(0)$ is the closure of the unit disc. The total dynamics of the agent are given by
\begin{equation}\label{eq:agent:dynamics}
    \begin{split}
        \dot{\agent}(t) &= f(\agent(t)) + u(t), \\
        \| u(t) \|_2 &\leq 1.
    \end{split}
\end{equation}

We consider a total of $M$ targets, each located at $\targetpos_i \in\ \interior\partition_{r_i}$ in some region $r_i$. We assume that each region contains at most one target, i.e., $r_{i_1} \neq r_{i_2}$ for $i_1 \neq i_2 \in \T = \{1,2,\dots, M\}$. Each target contains an internal state $\internal_i(t) \in \R^{m_i}$ following the stochastic \ac{lti} dynamics
\begin{equation}
    \dot{\internal}_i(t) = A_i \internal_i(t) + \omega_i(t),
\end{equation}
where $\omega_i$ is a zero-mean white noise process with covariance matrix $\internalvar_i \succ 0$. The interaction between the agent and target $i$ is described by the sensor measurement model
\begin{equation}\label{eq:outputs}
    \measurement_i(\agent(t), \internal_i(t)) = \monitoringfun_i(\agent(t)) H_i \internal_i(t) + \nu_i(t),
\end{equation}
where $\nu_i$ is a zero mean white noise process with covariance matrix $\measurevar_i$, and $\monitoringfun_i : \R^2 \to [0,1]$ captures the quality of the measurement depending on the position of the agent. While the precise form of $\gamma_i$ is not crucial, we assume it to (i) be differentiable with respect to the agent position within the region $\partition_i$ and (ii) have compact support $\support_i \subseteq \partition_i$.

For example, the internal states $\phi_i$ could describe a state of a natural disaster, e.g., the state of a flood or a fire. The outputs $\measurement_i$ would then be sensor measurements that allow us to estimate the considered state, e.g., water level at a dam, or the temperature of the fire at specific points of the target. The quality function $\gamma_i$ often has a dependence on the relative distance or angle between the agent and the monitored target, but its precise form strongly depends on the utilized sensors.

From well-known results we have that the optimal unbiased estimator of $\internal_i$ using the measurements from~\eqref{eq:outputs} is the Kalman-Bucy filter~\cite{lan2014variational}. The covariance matrix of this estimator is given by
\begin{equation*}
    \begin{split}
        \dot{\estimatecov}_i(t) &= A_i \estimatecov_i(t) + \estimatecov_i(t) A_i\super{\top} + \internalvar_i - \estimatecov_i(t) G_i(t) \estimatecov_i(t) \\ &\defeq f_{\estimatecov_{i}}(t)
    \end{split}
\end{equation*}
where $G_i(t) = \monitoringfun_i\super{2}(\agent(t)) H_i\super{\top} R_i\super{-1} H_i$. Although the observed system is in stochastic nature, the covariance matrices of the estimates are clearly deterministic and their evolution only depend on the agent position and the initial covariance matrices.

We assume that the pair $(A_i, H_i)$ is observable and $Q_i, \estimatecov_i(0) \succ 0$ for all $i \in \T$. These assumptions ensure that the covariance matrix $\estimatecov_i(t)$ is positive definite for all $t \geq 0$, and, that the covariance matrix trajectories eventually become periodic under any periodic control law that visits every target at least once~\cite{pinto2022multi}. We are interested in solving an infinite horizon problem; however, it is well known that periodic trajectories have the capability of approximating the optimal cost of average infinite horizon problems arbitrarily well~\cite{zhao2014optimal}. We thus restrict our attention to this subclass of solutions, and address the \ac{ocp} of minimizing the periodic average estimation error
\begin{mini!}
    {u, T}{\frac{1}{T}\int_0^T \sum_{i=1}^M \trace(\estimatecov_i(t))~\dt}{\label{ocp:global}}{\label{eq:average:loop:cost}}
    \addConstraint{\agent(0)}{= \agent(T)}
    \addConstraint{\estimatecov(0)}{= \estimatecov(T)}
    \addConstraint{\dot{\agent}(t)}{= f(\agent(t)) + u(t) \label{ocp:global:agent}}
    \addConstraint{\dot{\estimatecov}_i(t)}{= f_{\estimatecov_i}(t) \fat i}{\label{ocp:global:estimatecov}}
    \addConstraint{\| u(t) \|_2}{\leq 1,}
\end{mini!}
where the decision variable $T$ denotes the period. In the following section we characterize parameters that describe any optimal periodic trajectories.

\section{Trajectory Decomposition}\label{sec:decomposition}
Let $\agent^\ast : [0, T] \to \S$ be a periodic agent trajectory optimizing the cost function~\eqref{eq:average:loop:cost}. Then, $\agent^\ast$ is decomposed as follows:
\begin{enumerate}
    \item\label{item:decomposition:sequence} the sequence in which the targets are visited,
    \item\label{item:decomposition:trajectory:monitor} the \textit{monitoring trajectories} in each target region,
    \item\label{item:decomposition:trajectory:switch} the \textit{switching trajectories} connecting the end points of the monitoring trajectories with the respective start points of the subsequent trajectories,
    \item\label{item:decomposition:monitoring:duration} the \textit{monitoring durations} specifying the time spent in each target region, and
    \item\label{item:decomposition:switches} the \textit{switching points} given by the boundary points of the monitoring and switching trajectories.
\end{enumerate}
This decomposition, as illustrated in \figref{fig:decomposition}, allows rewriting the global cost function~\eqref{eq:average:loop:cost} in \textit{local} terms
\begin{equation}\label{eq:decomposition:global:cost}
    J = \frac{1}{T(\tau)}\sum_{k=1}\super{K} \bigl( \monitoringCT_k(\localT_k) + \integrationCT_k \bigr),
\end{equation}
where $\tau_k$ denotes the duration spent in the $k$th target region, $\monitoringCT_k$ denotes the local cost of the trajectory during the $k$th target visit (formally defined in \eqref{min:ocp:local:monitoring} below), $\integrationCT_k$ denotes the local cost of the trajectory switching between the $k$th and $(k+1)$st target visit (defined in~\eqref{eq:local:cost:switching} below), and $K$ denotes the length of the target visiting sequence. Note that the decomposition is quite heterogeneous: the trajectory segments are parameters in infinite-dimensional spaces; the monitoring durations and switching points are continuous parameters of finite dimension; and the visiting sequence is a sequence of discrete parameters of a priori unknown length. Solving~\eqref{ocp:global} directly and optimizing all parameters simultaneously is thus a very challenging problem. Instead, we propose a hierarchical decomposition approach. The discrete parameters are decided through a high-level (offline) planner, whereas the continuous parameters are optimized through a lower-level (online) planner. Let us briefly introduce all components of the decomposition in more detail.

\subsubsection{The Target Visiting Sequence}\label{sec:visiting:sequence} We proceed by first generating a visiting sequence $(i_1, i_2, \dots, i_K)$ that visits each target at least once. We do this by abstracting the mission space to a graph, where each node represents a target and each directed edge represents a time-optimal trajectory from one target to another.
The main challenge in abstracting the mission space is the generation of time-optimal trajectories between each pair of targets. To do this efficiently, we introduce the \ac{rrbt} algorithm in \secref{sec:global:planning}, a sampling-based algorithm tailored to the setting at hand. Given a root node, this algorithm generates a tree of nodes that allows us to approximate time-optimal trajectories from anywhere within the mission space to that node (see~\figref{fig:rrbt:illustration}). We generate such a tree for each target and use it to find approximately time-optimal trajectories between each target pair.

Based on this graph abstraction, we seek a minimum time cycle that visits all targets, i.e., we solve a \ac{tsp} problem. The \ac{tsp} cycle visits every node exactly once, with the exception being the start/end node. Note that it is possible that an edge from the TSP solution connecting target $i$ and $j$ traverses the region of another target $\ell$. To capture such revisits, we expand the TSP cycle to a closed path and insert all revisits. In the example, the visiting sequence would be $(\dots,i,\ell,j,\dots)$. A sample TSP solution is depicted in \figref{fig:tsp:initial:optimal:trajectory}, which is $(1,4,3,2)$. However, the region of target $3$ is revisited, which leads to a final visiting sequence $(1,4,3,2,3)$. To avoid confusion, we refer to the TSP solution as a \textit{cycle}, and to any closed trajectory that realizes the generated visiting sequence as a \textit{loop}. This entire task is performed off-line.

\subsubsection{The Monitoring Trajectories}
The monitoring trajectories describe the agent's path within each target region. For a given visit $i_k \in \T$ with \textit{entrance} and \textit{departure points} $\entrance_k, \departure_k \in \partial \partition_{i_k}$, monitoring duration $\localT_k$, and initial estimation covariance matrix $\bar{\Omega}$, we characterize the $k$th monitoring trajectory via the local \ac{ocp}
\begin{mini!}
    {\localT_k, \monitoringcontrol_k}{\int_{0}^{\localT_k} \sum_{i=1}^M \trace(\estimatecov_{i}(t))~\dt}{\label{min:ocp:local:monitoring}}{\monitoringCT_k(\localT_k) = }
    \addConstraint{\agent(0)}{ = \entrance_{k}}
    \addConstraint{\agent(\localT_k)}{ = \departure_{k}}
    \addConstraint{\estimatecov(0)}{= \bar{\estimatecov}}
    \addConstraint{\dot{\estimatecov}_i(t)}{= f_{\estimatecov_i}(t) \fat i \label{eq:ocp:local:monitoring:dynamics:estimate:cov}}
    \addConstraint{\dot{\agent}(t)}{ = f_{i_k}(\agent(t)) + \monitoringcontrol_k(t)  \label{eq:ocp:local:monitoring:dynamics:agent}}
    \addConstraint{\| \monitoringcontrol_k(t) \|_2}{\leq 1}
    \addConstraint{a(t)}{ \in \partition_{i_k}. \label{eq:ocp:local:monitoring:path:cosntraints}}
\end{mini!}
While analytical solutions may exist for some particular scenarios, we solve~\eqref{min:ocp:local:monitoring} using numerical optimization to achieve flexibility with respect to $f_{i_k}$ and the sensing function~\cite{rao2009survey}. This task is performed on-line.

\subsubsection{The Switching Trajectories}
The switching trajectories connect the switching points $\departure_k$ and $\entrance_{k+1}$ within the mission space. Since the agent does not sense any of the targets along these trajectories, and the covariance matrices are positive definite, the switching segments are time-optimal. The local cost of the switching period between the $k$th and $(k+1)$st monitoring period is
\begin{equation}\label{eq:local:cost:switching}
    \integrationCT_k = \int_{t^{\mathrm{sw}}_k}^{t^{\mathrm{sw}}_k
 + \Delta_k} \sum_{i=1}^M \trace(\estimatecov_i(t))~\dt,
\end{equation}
where $\Delta_k = \minswitchT(\departure_k, \entrance_{k+1})$ is the duration of the $k$th switching period, and $t^{\mathrm{sw}}_k$ denotes the starting time of the $k$th switching period. Although the agent may traverse multiple regions when transitioning from one target region to the next, none of these traversed regions contain a target. Thus, the dynamics of the estimator covariance matrices $\estimatecov_i$ are smooth during the switching segment.

\subsubsection{The Monitoring Durations}
The monitoring duration $\localT_k$ of the $k$th visit describes the amount of time the agent spends in the region $r_{i_k}$ of the $k$th target to collect measurements. In other words, it is the duration of the $k$th monitoring trajectory. To ensures that the local monitoring \ac{ocp}~\eqref{min:ocp:local:monitoring} is feasible, we require this parameter to be greater than the minimum time required to travel between the switching points $\entrance_k$ and $\departure_k$, i.e., we impose
\begin{equation}\label{eq:monitoring:duration:constraint}
    \localT_k \geq \minswitchT_{i_k}(\entrance_k, \departure_k).
\end{equation}
Those parameters are optimized in \secref{sec:optimization}. This task is performed on-line.

\subsubsection{The Switching Points}
As already introduced, we refer to the boundary points $\entrance$ and $\departure$ of the monitoring (and switching) trajectories as the switching points. In this paper, we fix these based on the TSP solution. Each edge of the TSP solution corresponds to a specific trajectory segment between two targets. We simply set the switching points to the intersection  of these trajectory segments with the boundary of the target regions. While it is conceivable to optimize those parameters in a similar fashion as the monitoring duration parameters (see~\secref{sec:optimization}), doing so is more involved, since it affects the switching trajectories and complicates the monitoring duration constraint~\eqref{eq:monitoring:duration:constraint}. Optimizing the switching points is left for future work.

In \secref{sec:optimization} we proceed with the optimization of the monitoring durations $\localT_k$ utilizing the decomposition into the local cost terms~\eqref{eq:decomposition:global:cost}. Before doing so, we introduce the \ac{rrbt} algorithm, which describes the missing piece of the higher level offline sequence planner.

\section{Global Planning}\label{sec:global:planning}
\begin{figure}
    \centering
    \includegraphics[width=\scalefactor\linewidth]{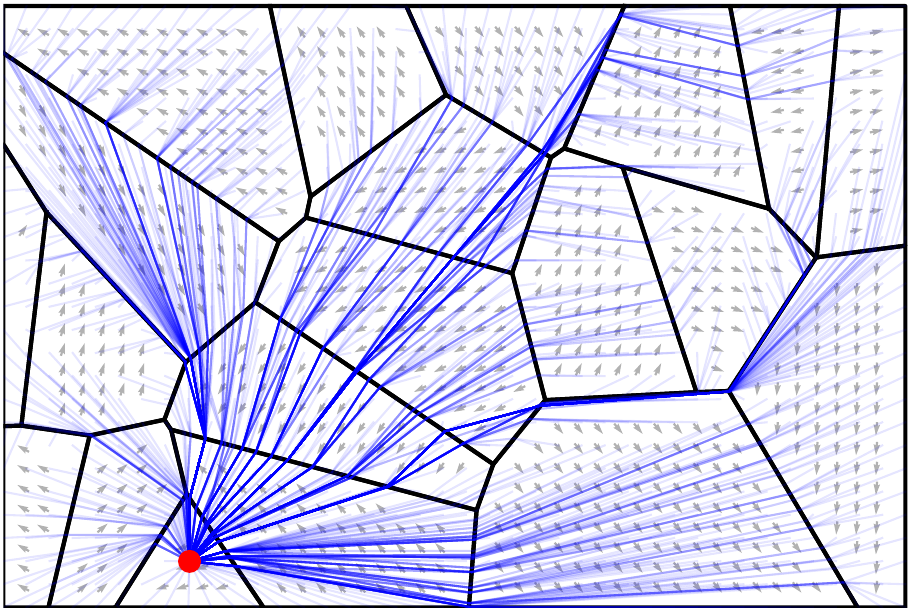}
    \caption{Illustration of global paths (blue) generated via the \ac{rrbt} algorithm for one root node (red).}
    \label{fig:rrbt:illustration}
\end{figure}

\begin{figure}
    \centering
    \includegraphics[width=\scalefactor\linewidth]{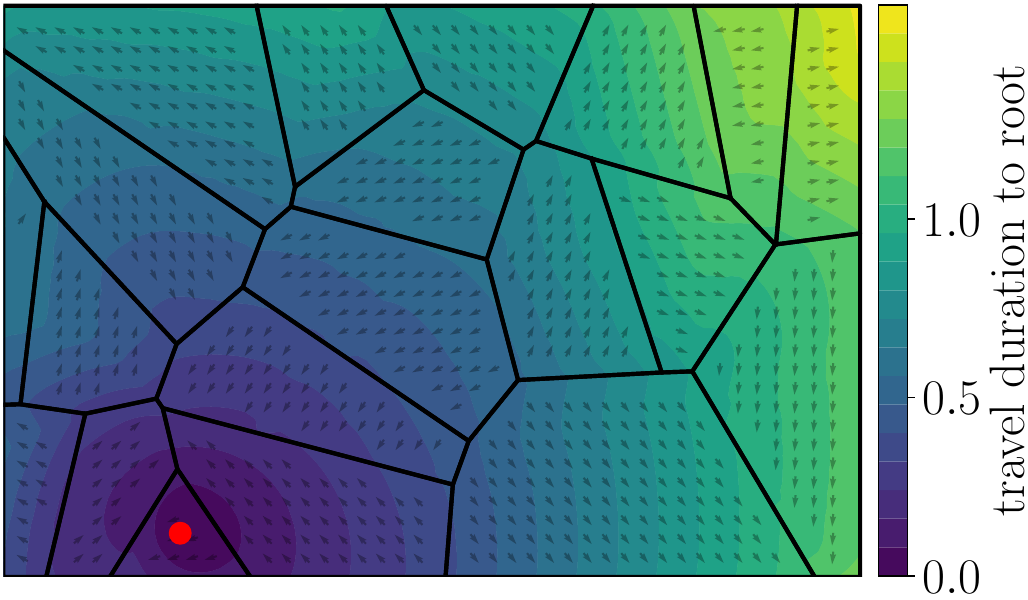}
    \caption{Illustration of travel duration from each location to the root node (red) using the \ac{rrbt} algorithm.}
    \label{fig:rrbt:illustration:travel:time}
\end{figure}

In order to generate the graph introduced in \secref{sec:visiting:sequence}, we must compute the minimum time required to travel between any pair of targets.
The main difference from the formulation in~\cite{karaman2010incremental} is the presence of the hybrid dynamics. We exploit the partition structure by assuming the existence of a local controller for each region $i$ with the capability of computing (time-optimal) control laws that steer the agent from $x$ to $y$ within region $\partition_i$. We denote the time required to do so by
\begin{equation}\label{eq:time:optimal:control:law}
    \minswitchT_i(x,y) : \partition_i \times \partition_i \to \R_{\geq 0} \cup \{ \infty \}.
\end{equation}
Note that while we assume $\minswitchT$ to be time optimal, the rest of our approach does not strictly require it; a heuristic could be used, with (potentially) a loss of optimality. With this abstraction, we can modify the well-known \ac{rrt} algorithm to produce a search tree on which all nodes, apart from the root node, lie on the boundary sets of the partition. Below we describe the \ac{rrbt} algorithm in more detail. A brief pseudocode overview is provided in \algoref{algo:rrbt}. \figref{fig:rrbt:illustration} illustrates an example of a search tree generated for the root node (red), and \figref{fig:rrbt:illustration:travel:time} shows the global travel durations from anywhere in $\S$ to the root node.

Given a goal position $x_f$ in the mission space, we initialize the tree of the \ac{rrbt} algorithm by creating a root node at $x_f$, and we initially \emph{activate} all regions which contain the root node. Active regions are all regions that have a representative within the tree; this is a necessary condition that guarantees a newly sampled node from that region can be connected to the tree by using the respective local controller~\eqref{eq:time:optimal:control:law}. Note that only one region is initially activated if the root lies in the interior of a partition set. Otherwise, multiple regions are activated, since each set of the partition is closed~\eqref{eq:partition:set}.

Let us now walk through a main iterate. The sampling phase aims to generate a new node to be included into the tree. We do this by first (uniformly) sampling a region $r$ from among the active regions. We then sample a boundary segment, and ultimately a new node $b$ from that segment. The neighborhood $N(b)$ of a boundary point $b$ is the set of all nodes in the tree that share a region with the sampled boundary point. Note that $N(b) \neq \emptyset$ by design. Given a potential parent $p \in N(b)$, we can compute an upper bound for the cost-to-root $c(b, x)$ by passing through $p$, i.e., $c(b,x) \leq \minswitchT_r(b, p) + c(p, x)$. We then connect the newly sampled node to a parent $\agent^\ast$ that minimizes this cost-to-root metric. The node is dropped if $\minswitchT_r(b,\agent^\ast) = \infty$. This procedure is repeated for a predefined number of iterations.  Once we have generated a search tree from a point $x$, we can find a global trajectory from any $y$ that lies in an active region by performing the connection phase of \algoref{algo:rrbt} with $b = y$. We use this mechanism to generate a tree for each target and then compute the target-to-target distance for each combination.

\def\notinbi{\bar{e}^i}
\def\notinbij{\notinbi_j}
\def\inbi{e^i}
\def\inbij{\inbi_j}
\def\neverinbi{\bar{E}^i}
\def\notneverinbi{E^i}
\def\notneverinb{E}
\def\rrbt{\rho}
\def\rrbtk{\rrbt^k}

The \ac{rrbt} algorithm converges in probability to a time-optimal path under mild conditions. To formalize this statement, let us introduce some notation. Although the index $k$ was previously used for the visiting sequence, let us use it in this section to denote the number of \ac{rrbt} iterations. We denote the duration of the best generated path from $x_0$ after $k$ iterates by $\rrbtk(x_0) \geq 0$, and the time-optimal travel duration by $\delta^\ast(x_0)$.
\begin{theorem}\label{thm:probabilistic:optimality:convergence}
    Let us assume that
    \begin{enumerate}
        \item\label{assumption:prob:completeness:finite:number:switching:points} there exists a time-optimal path $\agent^\ast : [0, \delta^\ast(x_0)] \to \R^2$ following a sequence of regions $(r_0, r_1,\dots,r_n)$ and boundary points $(s_1, s_2, \dots, s_{n})$, such that no boundary point borders more than three regions.
        \item\label{assumption:prob:completeness:time:optimal} the local switching durations $\minswitchT_i$ are time-optimal
        \item\label{assumption:prob:completeness:lipschitz} there exist neighborhoods $\N_i$ of $s_i$ for all $i=1,2,\dots,n$ such that $\delta_0\vert_{\{x_0\} \times \N_1}$, $\delta_n\vert_{\N_{n} \times \{x_f\}}$, and $\delta_i\vert_{\N_{i-1} \times \N_i}$ for $i=1,2,\dots,n-1$ are all Lipschitz continuous with Lipschitz constant $L$. Here, $\delta|_S$ denotes the restriction of $\delta$ to the set $S$.
    \end{enumerate}
    Then, for any $\varepsilon > 0$, it holds
    \begin{equation}\label{eq:probabilistic:optimality:convergence}
        \lim_{k \to \infty} \P \left(\rrbtk(x_0) - \delta^\ast(x_0) > \varepsilon \right) = 0.
    \end{equation}
\end{theorem}
A proof can be found in \appendixref{sec:appendix}.

\begin{figure}
    \centering
    \includegraphics[width=0.5\linewidth]{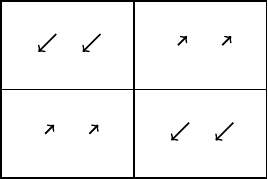}
    \caption{Illustrating the requirement of~\assumptionsref{assumption:prob:completeness:finite:number:switching:points}.}
    \label{fig:too:many:neighbors}
\end{figure}

\begin{remark}[Interpretation of~\assumptionsref{assumption:prob:completeness:finite:number:switching:points}]
    This assumption excludes scenarios where the optimal path switches between two regions that only have a single point in common. \figref{fig:too:many:neighbors} provides a counterexample in which the \ac{rrbt} algorithm will almost surely not find a feasible path. Assume that the vector fields in the second and fourth quadrant have large magnitudes. Then, any path connecting the third quadrant to the first quadrant must pass through the boundary point neighboring all regions, which has zero sampling probability.
\end{remark}

\begin{remark}[Interpretation of~\assumptionsref{assumption:prob:completeness:lipschitz}]
    The Lipschitz continuity of the local switching durations guarantees a bounded optimality gap for each local segment. Instead of imposing Lipschitz continuity only on the restrictions to the neighborhoods of the switching points, we could just assume (local) Lipschitz continuity. However, this would be too restrictive, since it would exclude all scenarios where $\minswitchT_i(y,z) = \infty$ for some pair $y,z \in S_i$.
\end{remark}

\begin{algorithm}
    \caption{RRBT* (from root node $x_f$)}
    \label{algo:rrbt}

    \algocomment{initialization phase} \\
    initialize tree with root at $x_f$ \\
    activate all regions containing $x_f$\\[0.5em]

    \algocomment{main loop} \\
    \While{\textnormal{max iter not exceeded}}{
        \algocomment{sample phase} \\
        sample $r$ from all active regions\\
        sample boundary point $b$ from $r$\\[0.5em]

        \algocomment{connection phase} \\
        get neighborhood $N(b)$ \\
        find parent $p \in N(b)$ minimizing cost-to-root \\
        connect $b$ to $p$ \\[0.5em]

        \algocomment{activation phase} \\
        activate any new region containing $b$
    }
\end{algorithm}

\section{Optimizing the Average Steady-State Cost}\label{sec:optimization}

We establish the global bilevel optimization problem
\begin{mini!}
    {\localT \in \R^K}{\frac{1}{T(\localT)}\sum_{k=1}\super{K} \bigl( \monitoringCT_k(\localT_k) + \integrationCT_k \bigr)}{\label{min:bilevel:optimization}}{\label{eq:cost:decomposition}}
    \addConstraint{\localT_k}{\geq \minswitchT_{i_k}(\entrance_k, \departure_k)}{ \fort k = 1, \dots, K, \label{eq:cost:decomposition:constraint}}
\end{mini!}
which couples the local segments of the decomposition described in \secref{sec:decomposition}. Note that the constraints~\eqref{eq:cost:decomposition:constraint} ensure that each local monitoring problem~\eqref{min:ocp:local:monitoring} is feasible, since the departure point $\departure_{k}$ becomes reachable from the entry point $\entrance_{k}$ (recall that $\minswitchT_{i_k}(\entrance_k, \departure_k)$ is the (minimum) time required to travel from $\entrance_k$ to $\departure_k$ in region $i_k$).

The gradient of the cost function~\eqref{eq:cost:decomposition} with respect to the parameter $\localT_k$ is given by
\begin{equation}\label{eq:gradient:bilevel:objective}
    \frac{\frac{d \monitoringCT_k(\localT_k)}{d \localT_k} T(\localT) - \sum_{k=1}^{K} \bigl(\monitoringCT_k(\localT_k) + \integrationCT_k\bigr)}{T(\localT)^2}.
\end{equation}
Note that the gradient $\frac{d M_k^\ast}{d \localT_k}$ is the \emph{sensitivity} of the $k$th optimal monitoring cost with respect to the parameter $\localT_k$. Assuming that strong duality holds in 
the monitoring OCP \eqref{min:ocp:local:monitoring}, we can utilize the shadow price equation~\cite{boyd2004convex}
$\frac{d M_k^\ast}{d \localT_k} = -\lambda_k$, where $\lambda_k \in \R$ is the dual variable of the equality constraint that fixes the monitoring duration $\localT_k$. This allows us to compute the gradient of the global cost function with respect to the local monitoring durations. We then solve \eqref{min:bilevel:optimization} in an alternating fashion, where we fix the monitoring durations and then simulate a loop. Note that simulating a loop includes re-solving each monitoring \ac{ocp} with the updated initial covariance matrices. We utilize the gradient \eqref{eq:gradient:bilevel:objective} and update the monitoring durations based on a simple projected gradient method with diminishing step size. We compare two optimization methods:
\begin{enumerate}
    \item\label{item:opti:variant:1} Simulating to steady state before updating $\localT$.
    \item\label{item:opti:variant:2} Updating $\localT$ upon completing each loop of the visiting sequence.
\end{enumerate}

The first method is motivated from the fact that the computed gradient~\eqref{eq:gradient:bilevel:objective} describes the steepest ascending direction of the global cost under a given configuration consisting of the monitoring durations and the estimator covariance matrices at the beginning of the loop. If the loop is not in steady state, the computed gradient may be inaccurate. On the other hand, the initial configuration, and thus the initial steady-state covariance matrices may be far from optimal, such that a precise gradient may be unnecessary to obtain. As we decrease the amount of change in the configuration, the loop will eventually converge to a steady-state loop, where the gradient will become precise. This motivates the second variant. Although convergence in bilevel optimization problems is difficult to establish, the next section empirically supports that both variants work well in practice.

\section{Numerical Results}\label{sec:numerical:results}

We discretize the local monitoring \ac{ocp}~\eqref{min:ocp:local:monitoring} via a multiple shooting method using the RK4 integrator~\cite{rao2009survey}. We model and solve this in CasADi~\cite{andersson2019casadi} utilizing IPOPT~\cite{wachter2006implementation}. Simulations were conducted on a hardware featuring an Intel i5 processor running at 1.60GHz with 16GB of RAM. All regions have randomly generated constant dynamics with norm bounded by one. Problem set up details can be found in the repository \url{github.com/hallfjonas/hytoperm}.

\begin{figure*}
    \centering
    \begin{subfigure}{0.5\linewidth}
        \centering
        \includegraphics[width=\scalefactor\linewidth]{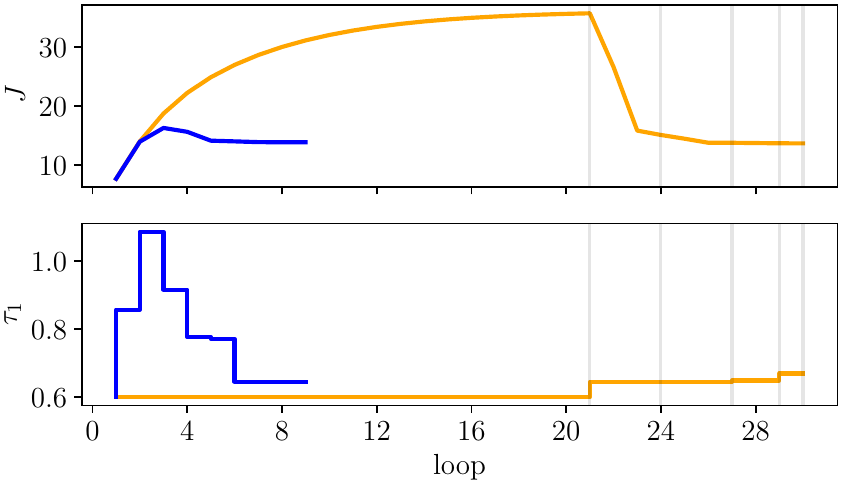}
        \caption{small experiment}
        \label{fig:optimization:small:experiment}
    \end{subfigure}%
    \begin{subfigure}{0.5\linewidth}
        \centering
        \includegraphics[width=\scalefactor\linewidth]{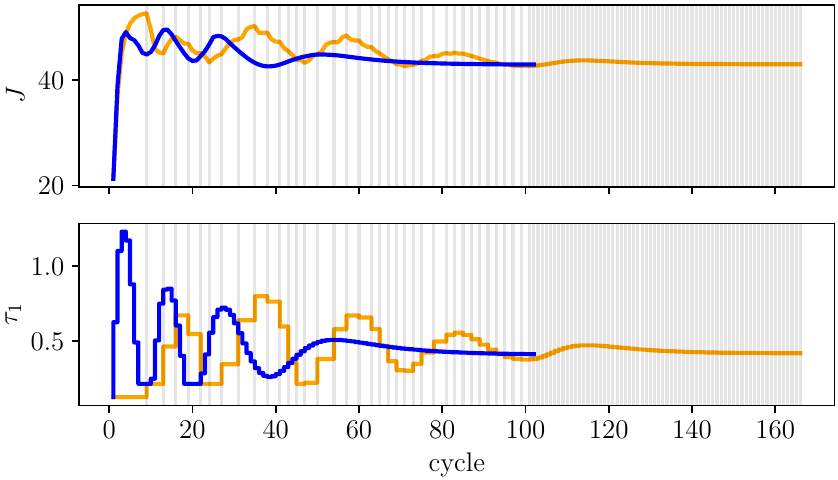}
        \caption{large experiment}
        \label{fig:optimization:large:experiment}
    \end{subfigure}
    \caption{Optimizing the monitoring durations $\localT$ utilizing both variant \ref{item:opti:variant:1} (yellow, vertical lines indicate steady-state) and \ref{item:opti:variant:2} (blue). The first plot shows the cost per loop, and the second plot shows the evolution of $\tau_1$.}
    \label{fig:optimization:experiments}
\end{figure*}

\begin{figure}
    \centering
    \includegraphics[width=\scalefactor\linewidth]{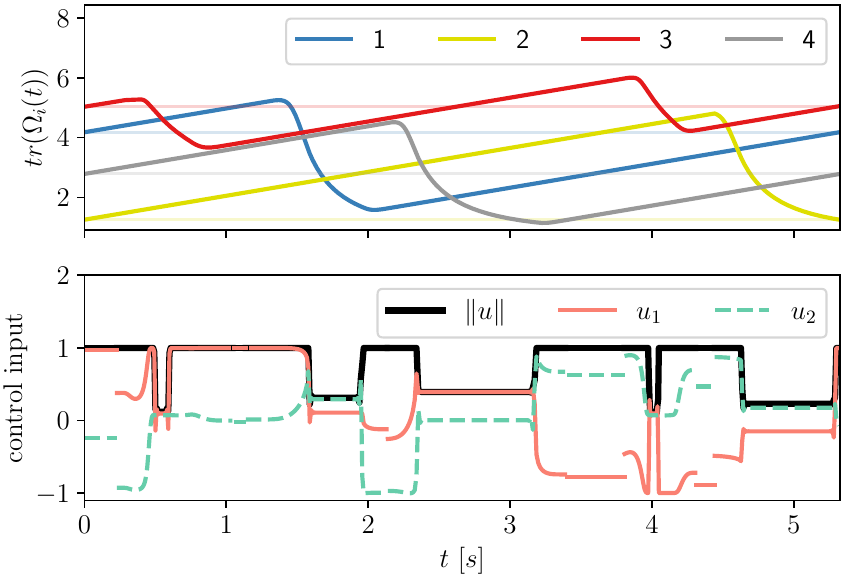}
    \caption{Depicting the optimized steady-state control law together with the mean estimation errors.}
    \label{fig:mse:and:controls}
\end{figure}

We first consider the scenario illustrated in \figref{fig:example:environment} with four targets and ten regions. \figref{fig:optimization:small:experiment} depicts the evolution of the global cost function together with the parameter $\tau_1$ to compare the two variants \ref{item:opti:variant:1} and \ref{item:opti:variant:2}. The evolution of the other monitoring durations have similar profiles. Note that in this particular scenario, variant \ref{item:opti:variant:2} completes its optimization procedure before variant \ref{item:opti:variant:1} even converges to the first steady-state loop. This is due to the fact that this variant is capable of modifying the monitoring periods more rapidly. As a result steady-state is reached more quickly. For both variants we see that the objective function is flat near the local optimum, since the cost does not change for the final few iterations while the monitoring durations are still updated. \figref{fig:mse:and:controls} depicts the controls for the optimal periodic trajectory together with the mean estimation errors for variant~\ref{item:opti:variant:1}. The results show that throughout the loop the control constraints $\| u(t) \|_2 \leq 1$ are respected, and that the mean estimation errors are indeed periodic. Both variants converge to the same loop (as depicted in~\figref{fig:tsp:initial:optimal:trajectory}).

A conceivable approach to improve the optimization procedure~\ref{item:opti:variant:1} is to implement an efficient method that solves the periodic Riccati equation and determines the steady-state trajectories of the covariance matrices given a fixed control law for a loop~\cite{varga2008solving}, rather than simulating the system to steady state. However, one could still not update the monitoring durations after each loop, since the dual variables of~\eqref{min:ocp:local:monitoring} would require updating as well. Those could be updated in a second loop. In total, this would lead to an optimization scheme where the monitoring durations are updated every other loop. For the scenario depicted in \figref{fig:optimization:small:experiment}, this would lead to a convergence within 10 loops (twice as many loops as the original method reaches steady state), which roughly compares to the number of loops required for variant \ref{item:opti:variant:2}. A detailed comparison is left for future work.

To show that the method extends to more complex scenarios, we consider another setting with 20 regions and 10 targets. \figref{fig:optimization:large:experiment} again demonstrates that directly updating the configuration after each loop reduces the number of loops required to reach the local optimum. Finally, we note that the cost iterations, even for variant \ref{item:opti:variant:1}, are not monotonic. This is due to the fact that the step size may be too large and result in a temporary cost increase. This is also reflected in the oscillation of the plotted monitoring duration.

\section{Conclusion and Future Work}\label{sec:future:work}
In this paper we presented a \ac{pm} formulation with the additional feature that the agent has piecewise smooth dynamics. With the goal of optimizing an average steady-state mean estimation error, we designed a periodic agent trajectory in two stages. First, in an offline step, we created a global path planner and computed a target visiting sequence. This yielded an initial closed trajectory for the agent to follow, which was optimized online.

As already mentioned, two aspects of the paper require a more detailed analysis. First, the switching segments are computed once and fixed afterwards. Their optimization in a global sense could benefit the overall performance. This is not a straightforward task, since even small alterations can change the traversed regions along the switching path. It is conceivable that suddenly another target region is visited along the altered switching path, which modifies the visiting sequence and requires the addition of a new monitoring segment.

\bibliographystyle{ieeetr}
\bibliography{bt/biblio}

\begin{appendices}
    \section{Proof of~\theoref{thm:probabilistic:optimality:convergence}}\label{sec:appendix}
    \begin{proof}
        Let us fix $\epsilon > 0$ and define $\beta = \frac{\varepsilon}{2(n+1)L} > 0$, where $n$ is the finite number of switching points and $L$ is the Lipschitz constant of the dynamics. Consider
        \begin{equation}\label{eq:boundary:neighborhood}
            B_i = B_{\beta}(s_i) \cap \N_i \cap \partial r_i \cap \partial r_{i-1},
        \end{equation}
        where $B_{\beta}(s_i)$ denotes the open ball of radius $\beta$ around $s_i$, and $\N_i$ is the neighborhood of $s_i$ in which the Lipschitz condition holds. By design, $s_i \in B_i$ and furthermore
        \begin{enumerate}
            \item\label{item:Bi:positive:measure} $B_i$ has positive measure in $r_i \cap r_{i-1}$, and
            \item\label{item:Bi:bounded:distance} $\| x - y \|_2 < \beta$ for all $x,y \in B_i$.
        \end{enumerate}
        We prove~\eqref{eq:probabilistic:optimality:convergence} in two main steps. First, in \claimref{claim:pc:zero:optimality:gap}, we show that if we sequentially sample from the sets $B_n, B_{n-1}, \dots, B_1$, then the \ac{rrbt} optimality gap diminishes. Then, in \claimref{claim:convergence:in:probability}, we show that the probability of sampling such a sequence converges to one. \claimref{claim:bounded:not:in:bi} can be seen as a lemma for the second claim.
        \begin{claim}\label{claim:pc:zero:optimality:gap}
            If there exists a sequence of iterates $i_n < i_{n-1} < \dots < i_1$ with $b_{i_\ell} \in B_\ell$, where $b_{i_\ell}$ is the sampled boundary point of the $i_k$th iterate of the \ac{rrbt} algorithm, then $\rrbtk(x_0) - \delta^\ast(x_0) \leq \varepsilon$ for all $k \geq i_1$.
        \end{claim}
        \begin{proof}
            Let us inductively show for $k=n,n-1,\dots,0$, that
            \begin{equation*}
                \rrbt^{i_k}(b_{i_k}) - \delta^\ast(s_k) \leq 2 (n-k+1) L \beta,
            \end{equation*}
            where we make the convention that $s_0 = b_{i_0} = x_0$. For the base case $k=n$, we have
            \begin{equation*}
                \rrbt^{i_n}(b_{i_n}) - \delta^\ast(s_n) = \delta_n(b_{i_n}, x_f) - \delta_n(s_n, x_f) \leq L\beta \leq 2L\beta.
            \end{equation*}
            At iterate $i_k$ we know that $b_{i_{k+1}}$ is already part of the tree. Thus, the \ac{rrbt} path from $b_{i_k}$ at iterate $i_k$ can be bounded by the path through the potential parent $b_{i_{k+1}}$. This provides
            \begin{equation*}
                \begin{split}
                    \rrbt^{i_k}(b_{i_k}) &\leq \delta_k(b_{i_k}, b_{i_{k+1}}) + \rrbt^{i_{k+1}}(b_{i_{k+1}}) \\
                    \delta^\ast(s_k) &= \delta_k(s_k, s_{k+1}) + \delta^\ast(s_{k+1}).
                \end{split}
            \end{equation*}
            Using Lipschitz continuity for both segments together with the induction step yields
            \begin{equation*}
                \rrbt^{i_k}(b_k) - \delta^\ast(s_k) \leq 2L\beta + 2(n-k)L\beta = 2(n-k+1)L\beta.
            \end{equation*}
            The claim follows when plugging in $\beta$ for $k=0$.
        \end{proof}

        Before proceeding, let us introduce some notation. Let $\notinbi_k$ denote the event that during the $k$th iterate region $r_i$ is sampled but $b_k \notin B_i$. Further, let $X_k$ be any event containing information about iterates prior to $k$, and let $\A_k$ denote the set of active regions at iterate $k$.
        \begin{claim}\label{claim:bounded:not:in:bi}
            If $\P(r_i \in \A_k \mid X_k) > 0$, then there exists a $\sigma < 1$ such that
            \begin{equation*}
                \P(\notinbi_k \mid r_i \in \A_k, X_k) \leq \sigma.
            \end{equation*}
        \end{claim}
        \begin{proof}
            If we condition on the number of active regions at iterate $k$ in addition to $r_i \in \A_k$, then the probability of $\notinbi_k$ is fully determined, and becomes independent of $X_k$. Let $\A_k^\ell$ denote the event that the number of active regions in iteration $k$ equals $\ell$, and define the total information vector $\Theta_k = \{r_i \in \A_k, X_k\}$.
            Then, by summing over all possible numbers of active regions given $\Theta_k$, we obtain
            \begin{equation*}
                \begin{split}
                    \P(\notinbi_k \mid \Theta_k) &= \sum_\ell \P(\notinbi_k \mid r_i \in \A_k, \Theta_k, \A_k^\ell) \P(\A_k^\ell \mid \Theta_k) \\
                    &= \sum_\ell \P(\notinbi_k \mid r_i \in \A_k, \A_k^\ell) \P(\A_k^\ell \mid \Theta_k) \\
                    &= \sum_\ell \frac{1-\alpha_i}{\ell} \P(\A_k^\ell \mid \Theta_k) \leq \frac{1-\alpha_i}{P} = \sigma < 1,
                \end{split}
            \end{equation*}
            where $\alpha_i > 0$ is the probability of sampling from $B_i$ given that region $r_i$ is sampled (positivity follows due to the finite number of boundary segments and~\ref{item:Bi:positive:measure}).
        \end{proof}

        We utilize this claim to address the question of how probable it is to sample such a sequence as described in~\claimref{claim:pc:zero:optimality:gap}. To this end, define the indices $I_i = (n-i+1)K$ for $K \in \mathbb{N}$. Then
        \begin{equation*}
            X^i_K = \cup_{k=I_{i+1}}^{I_i-1} \inbi_k
        \end{equation*}
        describes the event that \ac{rrbt} samples at least once from $B_i$ during the iterates $I_{i+1} \leq k < I_i$.
        \begin{claim}\label{claim:convergence:in:probability}
            For any $\gamma > 0$ there exists a $K$ such that
            \begin{equation*}
                \P\left( \cap_{i=1}^n X^i_K \right) \geq 1 - \gamma.
            \end{equation*}
        \end{claim}
        \begin{proof}
            We introduce $Y^i_K = \cap_{j=i}^{n} X^j_K$, the event of consecutively sampling at least once from $B_n$, then $B_{n-1}$, and so on, until eventually sampling from $B_i$ during the respective iterates. We note that each $Y^i_K$ only contains information about iterates prior to the index $I_i$. Further, it implies $r_{i-1} \in \A_k$ for all $k\geq I_i$, since $B_i \cap \partial r_{i-1} \neq \emptyset$. Repeated application of \claimref{claim:bounded:not:in:bi} yields $\P(X^i_K \mid Y^{i+1}_K) = 1 - \sigma^K$. Repeating this process for each $i$ provides
            \begin{equation*}
                \begin{split}
                    \P(Y^1_K) &= \P(X^1_K \mid Y^2_K) P(Y^2_K)= (1 - \sigma^K) P(Y^2_K) \\
                    &= \dots = (1 - \sigma^K)^n.
                \end{split}
            \end{equation*}
            Choosing $K > \log_\sigma(1 - \sqrt[n]{1 - \gamma})$ completes the claim.
        \end{proof}

        Now let $\gamma > 0$ arbitrary and $K$ such that $\P(\cap_{i=1}^n X^i_K) \geq 1 - \gamma$. Then, for any $k \geq I_1 = nK$ it holds
        \begin{equation*}
            \P \left(\delta_k - \delta^\ast < \varepsilon \right) \geq \P(\cap_{i=1}^n X^i_K) \geq 1 - \gamma.
        \end{equation*}
    \end{proof}
    \end{appendices}

\end{document}

%% file: ls/generic_article.tex
\input{ls/generic_setup.tex}
\input{ls/labels.tex}

\input{ls/theorems.tex}
\input{ls/acronyms.tex}

\usepackage{caption}
\usepackage{subcaption}

%% file: ls/generic_setup.tex
\usepackage{amssymb,amsfonts}
\usepackage{algorithmic}
\usepackage{mathtools}
\usepackage{textcomp}
\usepackage{hyperref}
\usepackage{lmodern}
\usepackage[utf8]{inputenc}
\usepackage[dvipdfm]{graphicx}
\usepackage{bmpsize}
\usepackage[export]{adjustbox}
\usepackage{cases}
\usepackage{empheq}
\usepackage{fancyhdr}
\usepackage{rotating}
\usepackage{bigdelim}
\usepackage{bm}
\usepackage{bbm}
\usepackage{tikz}
\usepackage{etoolbox}
\usepackage[short]{optidef}
\usepackage{multirow}
\usepackage{bookmark}
\usepackage{comment}

\newcommand{\secref}[1]{Sec.~\ref{#1}}
\newcommand{\appendixref}[1]{Appendix~\ref{#1}}

\newcommand{\theoref}[1]{Thm~\ref{#1}}
\newcommand{\claimref}[1]{Claim~\ref{#1}}
\renewcommand{\eqref}[1]{(\ref{#1})}
\newcommand{\figref}[1]{Fig.~\ref{#1}}

\newcommand{\algoref}[1]{Algorithm~\ref{#1}}

\newcommand{\assumptionsref}[1]{Assumption~\ref{#1}}

\newcommand{\myqed}{\hfill~\qed}

\newcommand{\R}{\mathbb{R}}

\newcommand{\defeq}{\vcentcolon=}

\newcommand{\ift}{\mathrm{~if~}}

\newcommand{\fort}{\mathrm{~for~}}

\newcommand{\fat}{\mathrm{~for~all~}}

\usepackage{upgreek}

\usepackage[ruled,vlined,linesnumbered]{algorithm2e}
\newcommand{\trace}{\mathrm{tr}}

\newcommand{\dt}{\textup{d}t}

\newcommand{\super}[1]{^{#1}}

\newcommand{\algocomment}[1]{{\textcolor{gray}{\# #1}}}

%% file: ls/labels.tex
\renewcommand{\theenumi}{(\roman{enumi})}%

%% file: ls/theorems.tex
\let\proof\relax 
\let\endproof\relax
\usepackage{amsthm}
\newtheorem{claim}{Claim}
\newtheorem{theorem}{Theorem}

\newtheorem*{assumptions*}{Assumptions}
\newtheorem{remark}{Remark}
\newtheorem*{task*}{Task}
\newtheorem*{problem*}{Problem}

%% file: ls/acronyms.tex
\usepackage[acronym,hyperfirst=false]{glossaries}

\newacronym[plural=LCQPs,firstplural=Linear Complementarity Quadratic Programs (LCQP)]{lcqp}{LCQP}{Linear Complementarity Quadratic Program}
\newacronym[plural=MPCCs,firstplural=Mathematical Programs with Complemenetarity Constraints (MPCC)]{mpcc}{MPCC}{Mathematical Program with Complemenetarity Constraints}
\newacronym[plural=QPs,firstplural=Quadratic Programs (QP)]{qp}{QP}{Quadratic Program}
\newacronym[plural=QCQPs,firstplural=Quadratically Constrained Quadratic Programs (QPQP)]{qcqp}{QCQP}{Quadratically Constrained Quadratic Program}
\newacronym[plural=NLPs,firstplural=Nonlinear Programs (NLP)]{nlp}{NLP}{Nonlinear Program}
\newacronym[plural=OCPs,firstplural=Optimal Control Problems (OCP)]{ocp}{OCP}{Optimal Control Problem}
\newacronym[plural=MIQPs,firstplural=Mixed Integer Quadratic Programs (MIQP)]{miqp}{MIQP}{Mixed Integer Quadratic Program}
\newacronym[plural=MINLPs,firstplural=Mixed Integer Non-Linear Programs (MINLP)]{minlp}{MINLP}{Mixed Integer Non-Linear Program}
\newacronym[plural=QPCCs,firstplural=Quadratic Programs with Complementarity Constraints (QPCC)]{qpcc}{QPCC}{Quadratic Program with Complementarity Constraints}
\newacronym[plural=QPLCCs,firstplural=Quadratic Programs with Linear Compementarity Constraints (QPLCC)]{qplcc}{QPLCC}{Quadratic Program with Linear Compementarity Constraints}
\newacronym[plural=FDIs,firstplural=Filippov Differential Inclusions (FDI)]{fdi}{FDI}{Filippov Differential Inclusion}
\newacronym[plural=ODEs,firstplural=Ordinary Differential Equations (ODE)]{ode}{ODE}{Ordinary Differential Equation}
\newacronym[plural=TSPs,firstplural=Traveling Salesperson Problems (TSP)]{tsp}{TSP}{Traveling Salesperson Problem}
\newacronym[plural=VRPs,firstplural=Vehicle Routing Problems (VRP)]{vrp}{VRP}{Vehicle Routing Problem}
\newacronym[plural=POIs,firstplural=Points of Interest (POI)]{poi}{POI}{Point of Interest}
\newacronym[plural=CLFs,firstplural=Control Lyapunov Functions (CLF)]{clf}{CLF}{Control Lyapunov Function}
\newacronym[plural=CBFs,firstplural=Control Barrier Functions (CBF)]{cbf}{CBF}{Control Barrier Function}
\newacronym[plural=OCBFs,firstplural=Optimal Control Barrier Functions (CBF)]{ocbf}{OCBF}{Optimal Control Barrier Function}

\newacronym{sqp}{SQP}{Sequential Quadratic Programming}
\newacronym{scp}{SCP}{Sequential Convex Programming}
\newacronym{licq}{LICQ}{Linear Independence Constraint Qualification}
\newacronym{mfcq}{MFCQ}{Mangasarian-Fromovitz Constraint Qualification}
\newacronym{kkt}{KKT}{Karush-Kuhn-Tucker}
\newacronym{pm}{PM}{Persistent Monitoring}
\newacronym{pmp}{PMP}{Pontryagin Minimum Principle}
\newacronym{lls}{LLS}{Linear Least Squares}
\newacronym{ipa}{IPA}{Infinitesimal Perturbation Analysis}
\newacronym{ad}{AD}{Automatic Differentiation}
\newacronym{mpc}{MPC}{Model Predictive Control}
\newacronym{rhc}{RHC}{Receding Horizon Control}
\newacronym{lti}{LTI}{Linear Time-Invariant}
\newacronym{rrt}{RRT*}{Rapidly-exploring Random Tree}
\newacronym{rrbt}{RRBT*}{Rapidly-exploring Random Boundary Tree}

\newcommand{\ac}[1]{\gls*{#1}}

%% file: ls/pm_specific_setup.tex
 \usepackage{bbm}

\renewcommand{\S}{\mathcal{S}}

\newcommand{\T}{\mathcal{T}}
\newcommand{\A}{\mathcal{A}}

\def\N{\mathcal{N}}

\def\P{\mathbb{P}}

\usepackage{booktabs}

\newcommand{\agent}{a}
\newcommand{\internal}{\phi}

\newcommand{\targetpos}{p}

\newcommand{\measurement}{z}

\newcommand{\estimatecov}{\Omega}
\newcommand{\internalvar}{Q}
\newcommand{\monitoringfun}{\gamma}

\newcommand{\measurevar}{R}
\newcommand{\partition}{S}

\newcommand{\partitionfunction}{g}

\newcommand{\entrance}{\agent\super{\varphi}}
\newcommand{\departure}{\agent\super{\psi}}
\newcommand{\localT}{\tau}

\newcommand{\integrationCT}{S^\ast}
\newcommand{\monitoringCT}{M^\ast}
\newcommand{\minswitchT}{\delta}

\newcommand{\monitoringcontrol}{u^{\mathrm{m}}}

\newcommand{\support}{K}
\newcommand{\interior}[1]{\stackrel{\circ}{#1}}

\usepackage{subcaption}

\newcounter{inlineenum}
\renewcommand{\theinlineenum}{\alph{inlineenum}}

\definecolor{agentcolor}{HTML}{326da8}
\definecolor{targetcolor}{HTML}{c84646}
